\documentclass[12pt]{amsart}
\usepackage{amssymb} 
\usepackage[all]{xy} 

\newcommand\ind{\operatorname{ind}}

\newcommand\datver[1]{\def\datverp% 
 {\par\boxed{\boxed{\text{Version: #1; Run: \today}}}}} 
\datver{1.1X; Revised: 03-11-98}

\newcommand\A{\mathcal A} 
\newcommand\E{\mathcal E}
\newcommand\B{\mathcal B}
\newcommand\Ah{{\mathbb A}^{\hbar}}  
\newcommand\CC{\mathbb C} 
\newcommand\fA{{\mathbb A}^ \hbar }
\def \n {\noindent}

\def \cal {\mathcal}
\def \uE { {}^{\E} }
\newcommand\fx{\hat{x}} 
\newcommand\fxi{\hat{\xi}} 
\newcommand\DD{\mathbb D}
\newcommand\NN{\mathbb N} 
 
\newcommand\RR{\mathbb R} 
\newcommand\ZZ{\mathbb Z} 
\def \cB {\overline{B}^*}

\newcommand\supp{\operatorname{supp}} 
\newcommand\Id{\operatorname{Id}} 
\newtheorem{theorem}{Theorem}

\newtheorem{proposition}{Proposition} 
 
\newtheorem{lemma}{Lemma} 
 
\theoremstyle{definition} 
\newtheorem{definition}{Definition} 
 
\theoremstyle{remark} 
\newtheorem{remark}{Remark} 
\begin{document} 
 
\title[Local formula for the  Index of FIO] 
{Local formula for the  Index of a Fourier Integral Operator} 

\author[E. Leichtnam]{Eric Leichtnam} \address{Institut de Jussieu-Chevelaret,
  Etage 7, plateau E,
  175 rue de Chevaleret 
        75013, Paris, France} 
        \email{leicht@math.jussieu.fr}
\author[R. Nest]{Ryszard Nest} 
        \address{Copenhagen}
       \email{rnest@math.ku.dk}
\author[B. Tsygan]{Boris Tsygan} 
        \address{Pennsylvania
        State University, Math. Dept., University Park, PA 16802, USA}
        \email{tsygan@math.psu.edu}

%\dedicatory\datverp 
%\date\datverp 
 
\begin{abstract}  
Let $X$ and $Y$ be two closed connected  Riemannian manifolds of the
same  dimension and  $\phi: S^*X \mapsto S^*Y$ a contact diffeomorphism.
  We show that 
the  index of an elliptic Fourier operator $\Phi$ associated with $\phi$
  is given by $\int_{B^*(X)}
  {\rm e}^{ \theta_0}\hat{A}(T^*X) 
 - \int_{B^*(Y)} {\rm e}^{ \theta_0}\hat{A}(T^*Y)$ where $\theta_0$ is
  a certain characteristic class depending on the principal 
symbol of $\Phi$  and $B^*(X)$ and $B^*(Y)$ are the unit ball bundles of
the manifolds  $X$ and  $Y$.  
 The proof uses the algebraic 
index theorem of Nest-Tsygan for symplectic  Lie Algebroids 
and an idea of Paul  Bressler to express the index of $\phi$ as a trace of $1$
in an appropriate deformed algebra. 

In the special case when $X=Y$ we obtain
 a different proof of  a theorem of Epstein-Melrose conjectured by Atiyah and 
 Weinstein.
\end{abstract}

\maketitle \tableofcontents 
 
\section{Introduction}

Let $X$ and $Y$ be two smooth closed connected Riemannian manifolds of the
same  
dimension  such that there exists a contact diffeomorphism 
$\phi: S^*X \mapsto S^*Y$ between the two unit cotangent bundles which
 induces a homogeneous symplectomorphism, still denoted by $\phi$, 
from $T^*X \setminus X$ onto $T^*Y \setminus Y$.

 We first recall the definition of the index of $\phi$ when 
dim$\,X \geq 3$  following
\cite{W 1}. We will denote by $\Omega_{{1 \over 2}} $ the half-density bundle over
$X$ or $Y$. Let $C_\phi$ be the graph of $\phi^{-1}$ in $(T^*Y\setminus
Y) \times (T^*X\setminus X)$ and $L_{C_\phi}$ be the associated Maslov
bundle. Let $A: L^2(X,\Omega_{{1 \over 2}}) \rightarrow L^2(Y,\Omega_{{1
\over 2}})$ be an elliptic Fourier Integral Operator of order zero whose
canonical relation is $C_\phi$ and whose principal symbol is an
invertible section of the bundle $\Omega_{1\over 2} \otimes L_{C_\phi}
\rightarrow C_\phi$ (see \cite{W 2}, \cite{Hormander}) for
details). Suppose that $B: L^2(Y,\Omega_{{1 \over 2}}) \rightarrow L^2(X,\Omega_{{1\over
2}})$ is an  elliptic Fourier Integral Operator of order zero   whose
canonical relation is $C_{{\phi}^{-1}}$. Then $B\circ A:\,
L^2(X,\Omega_{{1 \over 2}}) \rightarrow L^2(X,\Omega_{{1 \over 2}})$ is
an elliptic {\it scalar } pseudo-differential operator of order zero.
Since dim $X \geq 3$ there exists a smooth non vanishing 
function $x \in X \rightarrow a(x) \in {\CC}^*$ such that 
 the principal symbol of $ B\circ A$ is homotopic
to $(x,\xi) \in T^*X \rightarrow a(x) \in {\CC}^*$. In particular 
 the index of $ B\circ A$ is zero. Thus $ {\rm Ind} B = -
{\rm Ind} A $ for any Fourier Integral Operators $A$ and $B$ as above,
and, as the corollary of this fact, $ {\rm Ind} A$  does not depend on
the choice of $A$.  Since it only depends on the transformation $\phi$, it is called 
the index of $\phi$ and denoted by $ {\rm Ind} \, \phi$. A. Weinstein
has proved (see \cite{W 2}) 
that the integer $ {\rm Ind} \, \phi$ naturally appears if one wants to
compare the spectrum $(\,\lambda_k(X) \,)_{k \in \NN} $ of the Laplace
Beltrami operator $\Delta_X$ of $X$ with the one of $\Delta_Y$; for
instance if $T^*X \setminus X$ is simply connected then the sequence 
$
(\lambda_k(X) - \lambda_{k - {\rm Ind}\phi}(Y))_{k \in \NN}
$
is bounded.
 
The goal of this paper is to provide a geometric formula for the index an
elliptic  Fourier Integral Operator $\Phi$ of order zero 
 whose canonical relation is 
$C_\phi$ ( we do not assume dim$\,X \geq 3$).

\medskip

Let us first fix some notation. Given a smooth manifold $X$, we will use
$T^*X$ to denote the cotangent bundle of $X$ and $\cB X$ to denote the
projective compactification of $T^*X$. We will use $M$ to denote the  the
smooth manifold obtained by glueing at infinity  $ \cB X$ and $ \cB Y$  with
the help of the map
$$
\phi': (x, \xi) \rightarrow \phi(x,-\xi).
$$
Let $S^{0}(T^*X)$ and
  $S^{0}(T^*Y)$ denote the algebras of asymptotic symbols of
  pseudodifferential operators of order at most zero on $X$ and
  $Y$. Given an element $a\in S^{0}(T^*X)$, we denote by $a_\hbar$ the
   symbol  $a$ scaled by $\hbar$ in the cotangent
  direction and by Op($a$) the pseudodifferential operator associated
  to $a$. Given a pseudodifferential operator $A$, we denote by
  $\sigma (A)$ its full symbol (for the precise definition see the
  next section).

\medskip

The general strategy is as follows. We interpret conjugation by $\Phi$ as an
isomorphism of 
the algebras of pseudodifferential operators on $X$ and $Y$. Translated into
terms of formal deformations of the cotangent bundles, this allowes us to
construct a formal deformation $\fA (M)$ of $C^\infty (M)$ which on $T^*X$ and
$T^*Y$ represents the 
calculus of differential operators, while on the common cosphere at infinity
represents the calculus of pseudodifferential operators. While the symplectic
structures on $T^* X$ and $T^* Y$ do not glue together (so there is in general
no almost complex structure on $B^* X \cup_{\phi^{'}} B^* Y$), there is a
(noncanonical) symplectic Lie algebroid structure $(\E ,[\cdot ,\cdot ]
,\omega )$ over 
$M$ and  $\fA (M)$ is a 
deformation associated to it in the sense of \cite{N-T 2}. The usual traces on
the algebras of smoothing operators on $X$ and $Y$ give rise to a trace
$\tau_{\mbox{can}}$ on $\fA (M)$ such that $\ind \Phi = \tau_{\mbox{can}} (1)$.
An application of the general algebraic index theorem from \cite{N-T 2} gives
the local formula for the index. 

\medskip

The content of the paper is given below.

\medskip

\noindent {\bf 1.}

 In the first section we recall the relation between the
  calculus of $smoothing$  operators on $X$ and a formal
  deformation of $T^*X$ which is basically given by full symbol of a
  pseudodifferential operator.

\medskip

\noindent {\bf 2.}

${\cB}X$ carries a structure of symplectic Lie algebroid $(\E_X, [,], \omega)$
described in Section 2. The symbolic calculus of pseudodifferential operators
gives rise to a formal deformation of the sphere at infinity of ${\cB}X$
which, together with the formal deformation of $T^*X$ given above, gives rise
to a formal deformation $\fA (X)$ of ${\cB}X$ associated to $(\E_X, [,],
\omega)$.

\medskip

\noindent {\bf 3.}

Let us fix an almost unitary
 elliptic Fourier Integral Operator $\Phi$ whose canonical relation is
 given by the graph of $\phi^{-1}$.
In Section 4 we show how to glue together the deformations $\fA (X)$ and
$\fA (Y)$ into a formal deformation $\fA (M)$ of $M$ associated to a
symplectic Lie 
algebroid structure $(\E, [,], \omega)$ on $M$. The construction is based on
the following  
  strengthening of the Egorov theorem (see Theorem \ref{quantized}). 

\medskip

{\em \noindent 
\begin{enumerate} 
\item  
The map which to any 
   $a\in S^{0}(T^*X)$ associates the asymptotic expansion at $\hbar
   =0$ of 
 $ (\,\sigma ( \Phi {\rm Op}(a_\hbar) 
   \Phi^* )\,)_{{\hbar}^{-1}}$ induces an algebra 
   isomorphism 
$$
 \tilde{\Phi} :S^{0}(T^*X) \rightarrow S^{0}(T^*Y)
$$
\item  For each $k \in \NN^*$, there exists an $\E_X -$differential
operator $D_k$ on $\cB X $ such that, 
for any $a \in S^{0}(T^*X) $,  the following identity holds:
  $$
  \tilde{\Phi}( a  ) =( a + \sum_{k \geq 1} \hbar^k D_k(a ) ) \circ \phi^{-1} .
  $$ 
  \end{enumerate} }

\noindent Egorov theorem corresponds to the leading term in the
above expansion. 

The  real symplectic vector bundle $\E$  is isomorphic to $TM$ (as a vector
bundle over $M$) 
and hence  $TM$  is  the realification of a complex vector bundle
on $M$ which will be denoted by $\E_\CC$.

\medskip

\noindent {\bf 4}.

In Section 5 we identify the space traces on $\fA (M)$ and relate it to the
traces on the algebras of smoothing and of pseudodifferential operators.

\medskip

\noindent {\bf 5}.

In Section 6 we identify  the index of the Fourier Integral Operator with the
trace of 1 in the formal deformation. 

The local
index formula for $ {\rm Ind} \, \Phi$ follows from the algebraic index
theorem of \cite{N-T 2}, the class  
$\theta_0$ being the coefficient of $\hbar^0$ in the characteristic 
class $\theta$ ( \cite{N-T 2}, \cite{Fedosov}) of the deformation.

The main result can be formulated as follows.

{\em 
\begin{enumerate}
\item
Let $\Phi$ be a Fourier Integral Operator and $\fA (M)$ the formal deformation
of $M$ associated to it as in Definition \ref{RR}. Then
$$
\ind \Phi\,=\, \int_M {\rm e}^{\theta_0} \hat{A}(M) ,
$$
where $\theta_0$ denotes 
the characteristic class of the deformation of the 
 Lie Algebroid (${\mathcal E}, [\,,\,], \omega $) given by $\fA (M)$. 
\item
Let $\nabla_X$ be a connection $\nabla_X$ on the tangent bundle $T( \cB X)$
and $\hat{A}(T^*X)$ an associated representative form of  the
$\hat{A}-$class of $\nabla_X$.
The symplectomorphism 
$\phi$ induces a connection $\phi_* (\nabla_X)$ on the tangent space of 
${\cB} Y \setminus B^*(Y)$. 
Let $\nabla_Y$ denote its  extension to a connection of $T({\B Y})$
and $\hat{A}(T^*Y)$  an associated representative differential form  of the
$\hat{A}$-class of $\nabla_Y$. 
Then  
$$
 \ind \Phi\,=\, \int_{B^*(X)} {\rm e}^{\theta_0} \hat{A}(T^*X) \,-\,
\int_{B^*(Y)} {\rm e}^{\theta_0} \hat{A}(T^*Y) 
$$
   \end{enumerate}}  

\medskip

\noindent {\bf 6}.

The computation of the characteristic class of the deformation is
given in the last section, where we simultaneously construct a
deformation of $M$ and the Fourier Integral Operator whose index is
given by the trace of the 1 in the deformed algebra. As the
starting point we give  a somewhat nonstandard definition of the
characteristic class of a formal deformation which is
more amenable to computations in the case of deformations associated
to (twisted) differential or pseudodifferential operators.
 
As a corollary we get the following result.

{\em 
\begin{enumerate} 
\item There exists an almost unitary Fourier
  integral 
operator $\Phi_0$  whose canonical relation is
$C_\phi$ and such that:
$$
\ind \Phi_0 \,=\, \int_M \hat{A}(M)e^{{1 \over
    2}c_1(\E_\CC)}
$$
\item
If the dimension of $M$ is at least three, then 
$$
\ind (\phi ) = \int_M \hat{A}(M)e^{{1 \over 2}c_1(\E_\CC)}
$$ 
\end{enumerate}
 }
\noindent (compare with  \cite{E-M}, \cite{W 1}).

In the case when $X=Y$ a  straightforward Meyer-Vietoris type argument with
the mapping torus of $\phi$ shows that our results recover those of Epstein
and Melrose.

\medskip

\noindent {\bf 7}.

\begin{remark}
The methods of this paper extend in a fairly straightforward manner to the case
of a Fourier Integral Operator $\Phi$ between $L^2$ sections of vector bundles
$E$ and $F$
of the same dimension on $X$ and $Y$. In the case when both $X$ and $Y$ posses
a
metalinear structure, the corresponding index
formula is given by the expression 
$$
 \ind \Phi\,=\, \int_{B^*(X)} \mbox{ch}({\mathcal L}) \hat{A}(T^*X) \,-\,
\int_{B^*(Y)} \mbox{ch}({\mathcal L}) \hat{A}(T^*Y) 
$$
Here $\mathcal L$ is the vector bundle over $M$ obtained by glueing together
pull backs (by the canonical projections $\pi^*_X$ and $\pi^*_Y$) to the
cotangent bundles of 
$X$ (resp. $Y$) of the bundles 
${\Lambda^{\frac{n}{2}}}(X) \otimes E$ and ${\Lambda^{\frac{n}{2}}}(Y)
\otimes F$ with the help of the symbol of $\Phi$.

Note that existence of an
isomorphism of $\pi^*_X ({{\Lambda^{\frac{n}{2}}}}(X)
\otimes E)$ and $\pi^*_Y ({\Lambda^{\frac{n}{2}}}(Y)
\otimes F)$ over $\phi^{'}$ is equivalent to existence of an elliptic Fourier
Integral Operator from $L^2 (X, E)$ to $L^2 (Y, F)$.
\end{remark}
 
 \begin{remark}
 
It is easy to see that our local formula implies the following fact:

If $\phi$ extends as a symplectomorphism $:T^*X
\rightarrow T^*Y$ up to the zero section, then $ {\rm Ind} \,
\Phi\,=\,0$.  
  \end{remark}

\section{Symbolic calculus for $\Psi$DO's and formal deformations}

\subsection{Deformation of $T^*X$}

 \bigskip

 We will recall the pertinent facts from \cite{N-T 3}.

Let  $\chi$ be a smooth, non-negative  function on $X \times X$ satisfying the
following conditions: 

\medskip 
   (1) $\chi(x,y) = \chi (y,x)$;
  
   (2) $\chi \equiv 1$ on an open set containing the diagonal 
  in $X \times X$.
  
   (3) For each $x\in X$, the set $D_x=\{y\in X/\,(x,y) \in \supp \, \chi \}$
   is geodesically convex.
   
   \n We denote by ${\rm Exp}_x^{-1}$ the unique smooth inverse to the 
   exponential map:
   $$
   {\rm Exp}_x\,:\, T_x\,X \rightarrow X
   $$
   defined on $D_x$ and such that ${\rm Exp}_x^{-1}(x)=0$. 
   
   Given $x\in X,\, y \in D_x$, let $z$ denote the midpoint of the unique 
   geodesic joining $x$ and $y$ within $D_x$, and let $v\in T_z X$ be given by
   \begin{equation} \label{midpoint}
   v/2\,=\,{\rm Exp}_x^{-1} (y)\,=\,-{\rm Exp}_z^{-1} (x)
   \end{equation}
   Now, denote by $S^m(T^*X)$ the space of classical symbols of order $m$ on
   $X$, 
   i. e.. smooth functions $\theta$ on $T^*X$ satisfying estimates of the form:
   $$
   \sup_{(x,\xi)}\,|\partial_x^{\alpha}\,\partial^{\beta}_\xi\, \theta\,
   (x,\xi)|\, \leq \, C_{\alpha, \beta}\,(1 + |\xi|^2)^{{ m-|\beta| \over 2}}
   $$ 
$S^m(T^*X)$ is given the topology of (Frechet) topological vector space by the
 "best" $C_{\alpha, \beta}$.

 We will denote by
   $S^{+\infty}(T^*X) = \cup_{m\in \RR}S^m(T^*X) $ 
   the set of all classical symbols on $T^*X$.
   
   \n With the above notation (\ref{midpoint}), the map:
   $$
   {\rm Op}: S^m(T^*X) \rightarrow {\rm End}\,(C^{\infty}(X))
   $$
   given  by
   $$
   {\rm Op}(\theta)\,(u)\,(x)\, =\, \int_{T^*_zX}\, d\xi\, \int dy\, 
   \chi(x,y)\, {\rm e}^{i \,\xi.v}\, \theta(z,\xi)\, u(y)
$$
      defines a pseudo-differential operator. Conversely, if $P$ is a
   pseudo-differential operator on $X$  we define its complete 
   symbol to be:
   \begin{equation} \label{completesymbol}
   \sigma(P)(z,\xi)\,=\, P_y(\,\chi(x,y)\,{\rm e}^{i \,\xi.v} \,)|_{x=y=z}
   \end{equation}
   where $z$ is the midpoint of the geodesic joining $x$ and $y$ and 
$v$ satisfies (\ref{midpoint}). We observe that $P - {\rm Op}( \sigma(P)
)$ is a smoothing operator whose 
   Schwartz kernel vanishes to infinite
 order on the diagonal.
   Now, for
    a given $\theta \in C^{\infty}(T^*X)$, we  set:
   $\theta_\hbar(x,\xi)\,=\, \theta (x, \hbar \xi)$.

   Following (\cite{N-T 3}),
   we endow the algebra  
   $$\Ah (T^*X ) =
   C^{\infty}(T^*X) \otimes_\CC
   \CC[[\hbar]] $$ with a star product $\star_X$
    by defining,  for any symbols $\theta^1, \theta^2, \in S^m(T^*X)$,
   $\theta^1 \star_X \theta^2$ to be the asymptotic expansion at
   $\hbar=0$ of:
 \begin{equation} \label{star}
 \sigma ( \,{\rm Op}(\theta^1_\hbar) \circ
   {\rm Op}(\theta^2_\hbar)\,)_{\hbar^{-1}}
   \end{equation}

 One sees immediately that $\star_X$ extends to $\Ah (T^*X )$

 Recall 
   that there exists, unique up to normalization, a 
 canonical trace on $(\Ah (T^*X ),
   \star_X)$, ${\rm Tr}_{{\rm can}}^X$,  given  
   by:
   $$
   \forall a \in S^{-\infty}(T^*X),
   \quad {\rm Tr}_{{\rm can}}^X\,(a)\,=\,
    {\rm Tr} (\,{\rm Op} (a_\hbar) \,)\,=\, 
   { 1 \over n! \, \hbar^n}\, \int_{T^*X} a\, (\omega^X)^n\in 
   \CC[\hbar^{-1}, \hbar]]
   $$ ( Proposition 2.5 (3) of \cite{N-T 3}  ).

\subsection{Lie algebroid structure and deformation quantization the
  projective completion $\cB X$}

For any $x\in X$, we set $  {B}^*_x X = \frac{\RR_+ \oplus T^*_x X} {\RR_+^*}
 \setminus \{(0,0)\} $, and embed $T^*_x X$ in ${\cB }_x X$ by sending $\xi$
 to the class  
  of $1\oplus \xi$.  We view 
  $ {B}^*_x X$ as a compactification of $T^*_xX$. Then  we consider the
  fiber bundle $\cB X$ over $X$ defined by
   $\cB X= \cup_{x\in X}{B}^*_x X $. Therefore $\cB X$ is a compactification
   of $T^* X$ and a  smooth compact
   manifold with boundary: $\partial \cB X = \cB X \setminus T^*X$. 
   Similarly one defines the bundle $\cB Y$ over $Y$.  
  We observe 
   that the map from $S^*X$ into  $\cB X$ given by $\xi \rightarrow 0\oplus 
   \xi$ defines an isomorphism between $S^*X$ and $\cB X\setminus
 T^*X$. For any  $\xi=(x,\xi_x) \in T^*_xX$ we will define 
 $-\xi$ to be $(x,-\xi_x)\in T^*_xX $.
    Clearly,  
   $\phi$ induces a natural smooth isomorphism of manifolds 
   with boundary:
   $$\phi': \cB X \setminus X \mapsto \cB Y\setminus Y$$
 defined 
   by 
$$\phi' (\lambda \oplus \xi) = (\lambda \oplus \phi(-\xi) )\;
 {\rm { if}}\, \xi \in T^*X \setminus X,\; 
 \phi' (\lambda \oplus 0) = (\lambda
   \oplus 0 ).$$
 
 By
     glueing
   $\cB X$ and $\cB Y$ along the boundary $\cB X \setminus T^*X$ with the help
   of $\phi'$, we define the following  smooth compact manifold $M$: 
   \begin{equation} \label{glue}
   M = \cB X \cup_{\phi'} \cB Y
   \end{equation}
   
    Let $\Pi_X: \cB X \rightarrow X$ 
   be the projection map.  We denote by $\Xi^X$ the set of smooth 
    vectors fields of $\cB X$ which are tangent 
   to all the submanifolds $\Pi_X^{-1}(x) \cap (\cB X \setminus T^*X)$, $x\in
   X$. Let $(x,\xi)=(x_1,\ldots, x_n; \xi_1,\ldots, \xi_n)$ 
   be a local chart of $T^*X$ and $(\rho, \theta)$ be 
   the polar coordinates: $ \rho = || \xi ||, \;
   \theta = { \xi \over || \xi ||} $, where $ || \, ||$ denotes the Euclidean norm of $T^*X$. Then a local chart 
   of $\cB X$ near $\cB X \setminus T^*X$ is given by
   \begin{equation} \label{chart}
   (x_1,\ldots, x_n; t= { 1\over \rho}, 
   \theta=(\theta_1,\ldots, \theta_{n-1} ))\; t \geq 0,\; \theta \in S^{n-1}
   \end{equation}
   In this local chart, $\Xi^X$ is generated by the 
   vector fields $t{\partial \over \partial x_j}$, 
   $t{\partial \over \partial t}$, ${\partial \over \partial \theta_l} $, where
$1\leq j \leq n$, $1\leq l \leq n-1$. We will use several times the
following obvious Lemma
   \begin{lemma} \label{lemma}
The vector fields ${\partial \over \partial \xi_j} \, \, (1 \leq j \leq
n)$ belong 
to the $C^{\infty}(\cB X)-$module $t \Xi^X$ generated by $t^2 {\partial
\over \partial t}$, 
   $t{\partial \over \partial \theta_l} \, \,(1\leq l \leq n-1)$.
   \end{lemma}
   
    Moreover 
we observe that the set of classical symbols of order zero on $T^*X$ is
nothing else 
   but $C^{\infty}(\cB X)$. 
   
   Before we continue, let us recall the definition of a 
symplectic Lie algebroid  (
see for instance \cite{Mac Kenzie}, \cite{ N-T 2}). 
\begin{definition} {\item 1)} A symplectic Lie algebroid on $M$ is a quadruple 
$(\E, \rho,  [\, , \,], \omega)$ on $M$, where $\E$ is a smooth vector
bundle on $M$, $[\, , \,] $ is a Lie algebra structure
 on the sheaf of sections
of $\E$, $\rho$ is a smooth map of vector bundles:
$$
\rho: \E \rightarrow TM
$$ such that the induced map:
$$
\Gamma(\rho): C^{+\infty}(M,\E) \rightarrow  C^{+\infty}(M, TM)
$$ is a Lie algebra homomorphism and, for any sections $\sigma$ and 
$\tau$ of $\E$ and any smooth function $f$ on $M$, the following 
identity holds:
$$
[\, \sigma, f \tau\,]\,=\, \rho(\sigma)(f).\tau + f[ \sigma, \tau]
$$ Lastly, $\omega$ is a closed $\E-$two form on $M$ such that the
associated linear map:
$$
C^{+\infty}(M,\E) \times C^{+\infty}(M,\E) \ni (U,V)
 \mapsto \omega(U,V) \in C^{+\infty}(M)
$$ defines a symplectic structure on $\E$.
{\item 2)} The ring of $\E-$differential operators is by definition
the ring generated by smooth functions on $M$ and smooth  sections of $\E$.
{\item 3)} We denote by ${}^{\E}\Omega^{{\bf{.}}} \,=\, C^{+\infty}( M,
\Lambda^{{\bf{.}}}\E^*)$ the set  of smooth 
 sections on $M$ of the bundle of  alternating multilinear forms on $\E$.
\end{definition}
   
   We leave to the reader the  easy proof 
   of the following:
   \begin{proposition} \label{algebroid} For any $p \in \cB X$, we set 
   $$ 
   {\E}_{p}^X\,=\, { \Xi^X \over I_p \Xi^X}
   $$ where $I_p$ is the set of smooth real-valued functions on $\cB X$ which 
   vanish at $p$. Then:
   {\item 1)}  $({\E}_{p}^X)_{p\in \cB X}$ 
   form a smooth vector bundle, denoted $ {\E}^{X}$, over $\cB X$ such that 
   the set of smooth sections over $\cB X$ of $ {\E}^{X}$ 
is the same as
 $\Xi^X$. If $U,V \in \Xi^X$ then the 
Lie
bracket $[U,V]$ also 
   belongs to $\Xi^X$.
   {\item 2)} The fundamental two-form $\omega^X (= \sum_{j=1}^n\, d\xi_j \wedge 
   d x_j)$ of $T^*X$ induces a smooth form, still denoted $\omega^X$, in 
   $C^{\infty}( \cB X; {\uE}^{X}\Omega^2\,)$. Moreover, $({\E}^{X},[,],\omega^X)$
   defines a symplectic Lie algebroid   over $\cB X$. 
  \end{proposition}
   
   \begin{proposition} \label{ext} The star product $\star_X$ on $T^*X$ extends 
   to a star product, still denoted  $\star_X$, on 
   $\cB X$ such that for any $f,g \in C^{\infty}(\cB X )$ we have:
   $$
   f\star_X g = f g + \sum_{ n \geq 1} \hbar^n A^{(n)}(f\,,\,g)
   $$ where the $A^{(n)}$ are ${\E}^X-$bidifferential operators.
   \end{proposition} 
   \begin{proof} Let $(x,\xi)=(x_1,\ldots, x_n; \xi_1,\ldots, \xi_n)$ 
   be a local chart  of $T^*X$. Then for any $f,g \in C^{\infty}(\cB X )$ 
   and $(x,\xi)$ in the domain of this local chart we have:
   $$ 
   f\star_X g (x, \xi)=
 \sum_{\alpha, \beta \in {\NN}^n,\, |\beta|\leq |\alpha|} \, 
   {\hbar^{|\alpha|} \over \alpha !} \, c_{\alpha, \beta}(x)
   D_\xi^{\alpha}f(x, \xi)  {\partial^{\beta} 
   \over \partial^{\beta} x}g(x, \xi)
   $$
then, using the local coordinates (\ref{chart}) and Lemma \ref{lemma},
one gets
    easily all the results of the proposition.
     \end{proof}
     
      Proposition \ref{algebroid} allows to formulate the following
      definition.

    \begin{definition} \label{almostcomplex}
 {\item 1)} A  smooth real-vector
     bundle ${\E}^Y$ over 
    $\cB Y$ is defined by setting
    ${\E}^Y_{|{T^*Y} }= T(T^*Y)$ and 
    ${\E}^Y_{|{\cB Y \setminus Y}}= \phi_*({\E}^X_{|{\cB X \setminus
    X}})$.
     {\item 2)} By
 glueing  ${\E}^X$ and ${\E}^Y$ along $\cB X \setminus T^*X$
     with the help of $\phi'$, one defines a smooth vector bundle $\E$ 
     over $M$ which is isomorphic to $TM$. A  smooth exact  differential 
     form $\omega \in C^{\infty}( M; {\uE}\Omega^2\,)$ is defined by
     setting $\omega_{|\cB X} = \omega^X$, 
     $\omega_{|\cB Y \setminus Y} = 
     (\phi^{-1})^*(\omega^X_{|\cB X \setminus X})$ and 
     $\omega_{|T^*Y }=\omega^Y $ (where $\omega^Y$ is the 
canonical two form of $T^*Y$). 
 {\item 3)} The natural injection:
$$
C^{+\infty}(M, \E) \rightarrow C^{+\infty}(M, TM)
$$ is induced by a bundle map $\rho: \E \rightarrow TM$ as in 
Proposition \ref{algebroid} and
 $({\E},\rho, [,],\omega)$ 
   defines a symplectic Lie algebroid  which will be denoted 
 $({\E}, [,],\omega)$ in the sequel.
\end{definition}
%%%%%%%%%%%%%%%%%%%%%%%%%%%%%%
%%%%%%%%%%%%%%%%%%%%%%%%%%%%%
%%%%%%%%%%%%%%%%%%%%%%%%%%%%% 

\section{Regularized Index formula for a Fourier Integral Operator}

%%%%%%%%%%%%%%%%%%%%%%%%%%%
%%%%%%%%%%%%%%%%%%%%%%%%%%%
Let $C_\phi$  be
  the graph of $\phi^{-1}$ 
 in $(T^*Y\setminus Y) \times (T^*X\setminus X)$ and $L_{C_\phi}$ be the 
 associated Maslov bundle over $C_\phi $. 
 We fix $\Phi: L^2(X,\Omega_{{1 \over 2}}) \rightarrow L^2(Y,\Omega_{{1 \over
 2}})$ an elliptic Fourier integral 
 operator of order zero whose canonical relation is $C_\phi$   and
  whose principal 
 symbol $a$ is a unitary section   of the bundle 
$\Omega_{1\over 2} \otimes L_{C_\phi} \rightarrow C_\phi$: this means
that $a$ is homogeneous
 of degree zero (i.e. constant on each ray) and 
 that $a \overline{a} \equiv 1$: see \cite{W 2}.  {\it We can, and will,
  assume in the sequel that 
 $\Phi \Phi^* - \Id$ and $\Phi^* \Phi - \Id$ are smoothing}.
 As observed in \cite{W 2} 
 $\Phi$ is Fredholm, with index  defined by $\ind \Phi = 
 {\rm dim}\,\ker \Phi -
 {\rm dim}\, {\rm coker} \Phi$. In order to give a formula "via
 regularization" for $\ind \Phi$ we introduce the following algebra $\A$
 which will have a "regularized" trace.
 $$
 \A = \{ (A,B) \in \Psi^0(X, \Omega_{{1\over 2}} ) \times \Psi^0(Y,
 \Omega_{{1\over 2}} )|\; A -\Phi^* B \Phi\, {\rm is}\, {\rm smoothing}\}
 $$

 We leave to the reader the easy proof of the following:

 \begin{proposition} \label{regtrace} 
\item{ 1)} 
The map 
 $\tau: \A \rightarrow \CC$ given by:
 $$
  \forall (A,B) \in \A,\; \tau (A,B)\,
  =\, {\rm Tr} ( A -\Phi^* B \Phi) -  {\rm Tr}( B ( \Id - \Phi \Phi^*) )
  $$ 
is a trace
  \item{ 2)} 
$\ind \Phi = \tau( \Id, \Id)$.
 \end{proposition}
  
 \begin{remark} 
$\tau(A,B)$ is a "regularization" of 
${\rm Tr}  A - {\rm Tr}  B$.
 
 \end{remark}
  
 \section{Algebraization of a Fourier Integral Operator} 
 
 We are going to use the  following (deformed quantized algebra), where the 
 manifold 
 $Z$ is equal to $X$ or $Y$:
 $$
{\mathbb B}^\hbar (\cB Z ) = { C^{\infty}(\cB Z) \over 
 C_0^{\infty}(\cB Z)}\otimes_\CC \CC[[\hbar]]
 $$ 
 where $ C_0^{\infty}(\cB Z)$ denotes the set of smooth functions
 which vanish to infinite order at $ \cB Z \setminus T^*Z$. We observe 
 that  $\star_Z$ induces a star-product, still denoted $\star_Z$, on
 ${\mathbb B}^\hbar (\cB Z )$.
   
   \begin{theorem} \label{quantized} 
{\item 1)}  The map which to any 
   $a\in S^{0}(T^*X)$ associates the asymptotic expansion at $\hbar
   =0$ of 
 $ (\,\sigma ( \Phi {\rm Op}(a_\hbar) 
   \Phi^* )\,)_{{\hbar}^{-1}}$ induces an algebra 
   isomorphism $\tilde{\Phi}$
    from $({\mathbb B}^\hbar (\cB X ), \star_X)$ onto
     $({\mathbb B}^\hbar (\cB Y ), \star_Y)$.
{\item 2)} For each $k \in \NN^*$, there exists an $\E-$differential
operator $D_k$ on  $\cB X$ such that 
for any $a \in C^{+\infty}( \cB X)$ which is identically zero in a
neighborhood of the zero 
  section, we have the following identity :
  $$
  \tilde{\Phi}( a  ) =( a + \sum_{k \geq 1} \hbar^k D_k(a ) ) \circ \phi^{-1}
  $$  in the vector space ${\mathbb B}^\hbar (\cB Y )$  
  \end{theorem}

 Before proving this theorem we state the next
  proposition which is 
 an easy consequence of Proposition \ref{ext} 
   and Theorem \ref{quantized}
\begin{proposition} \label{star5} 
The star product $\star_Y$ on $T^*Y$
extends 
   to a star product, still denoted  $\star_Y$, on 
   $\cB Y$ such that for any $f,g \in C^{\infty}(\cB Y )$ we have:
   $$
   f\star_Y g = f g + \sum_{ n \geq 1} \hbar^n B^{(n)}(f\,,\,g)
   $$ where the $B^{(n)}$ are ${\uE}^Y-$bidifferential operators.
   \end{proposition}

  \begin{proof} Let us first assume
part 2). Then,  using the results of Section 2.1 and  the fact
that $\Phi \Phi^* - {\rm Id}$ and $ \Phi^* \Phi- {\rm Id}$ are
smoothing, one proves 
  easily that $\tilde{\Phi}$ is an isomorphism whose inverse is given by:
  $$
  b\in S^{0}(T^*Y) \rightarrow (\,\sigma ( \Phi^* {\rm Op}(b_\hbar) 
   \Phi )\,)_{{\hbar}^{-1}}
$$ Now let us prove part 2). Following \cite{Hormander} page 26, we
recall that the Schwartz kernel 
of $\Phi$ is the finite sum of a smooth function and of oscillatory
integrals (supported in small 
  coordinates charts) of the following type:
  \begin{equation} \label{Schwartz1}
 K(y,x)\,=\, \int_{\RR^n} e^{i ( \varphi(y, \eta) -x\cdot\eta)} \,b(y,\eta) d\eta
  \end{equation}
where $ b(y,\eta) \in S^0(T^*Y)$ vanishes for $|| \eta|| \leq 1$,
$\varphi(y, \eta)$ is
an homogeneous phase function parametrizing locally the graph $C_\phi$
of $\phi^{-1} $ which satisfies 
$ {\rm det} { \partial^2 \varphi \over \partial y \partial \eta } \not=
0$ so that locally we have:
   $$
\{ (y,\varphi'_y(y, \eta)\,;\,\varphi'_\eta (y, \eta), \eta
)\}\,=\,C_\phi 
$$ and $\phi^{-1}(y,\varphi'_y(y, \eta ) )= (\varphi'_\eta (y, \eta),
\eta )$. Notice moreover 
   that $(y, \eta) \rightarrow (y,\varphi'_y(y, \eta ) )$ and 
$(y, \eta) \rightarrow (\varphi'_\eta(y, \eta), \eta )$ are local
diffeomorphisms. 
 
 With these notations,   the Schwartz kernel 
of $\Phi^*$ is the finite sum of a smooth function and of oscillatory
integrals (supported in small 
  coordinates charts) of the following type:
  \begin{equation} \label{Schwartz2}
K^*(x,y)\,=\, \int_{\RR^n} e^{-i ( \varphi(y, \eta) -x\cdot\eta)}
\,b_1(y,\eta) d\eta
  \end{equation}
  
Let $a \in C^{+\infty}( \cB X)$ which is identically zero in a
neighborhood of the zero 
section, in order to analyze $\Phi \circ {\rm Op}(a_\hbar)$ it is enough
to study the 
operator $K \circ {\rm Op}(a_\hbar)$ where $K$ denotes the operator
whose Schwartz kernel is 
  given by (\ref{Schwartz1}). 
  
  The Schwartz kernel of $K \circ {\rm Op}(a_\hbar)$ is given by:
  $$
  T(y,z) = \int_{\RR^n} \int_{\RR^n} \int_{\RR^n} 
e^{i ( \varphi(y, \eta) -x\cdot\eta)} \,b(y,\eta) \, a(x, \hbar \xi) e^{i
(x-z )\cdot \xi } d\xi d x d\eta
  $$ 
In this integral we replace $ a(x, \hbar \xi)$ by its Taylor expansion: 
$$ 
\sum_{\alpha \in \NN^n} { 1 \over \alpha !} \partial^{\alpha}_x a(z,
\hbar \xi) (x -z)^{\alpha}
  $$ 
Using 
  the following two identities
 
  $$
(x -z)^{\alpha} e^{i (x-z )\cdot \xi }\,=\, D^{\alpha}_\xi (e^{i (x-z
  )\cdot \xi }
)
  $$ 
$$ \int_{\RR^n} e^{i x\cdot (\xi -\eta) } dx \,=\, (2\pi)^n
\delta_{\xi = \eta}
$$
and integrating by parts we see that $T(y,z) $ is the sum of a smooth
function and of

 $H(y,z)=$
  $$
   \int \int \int 
  e^{i ( \varphi(y, \eta) -x\cdot \eta)} \,b(y,\eta)  \,
   \sum_{\alpha \in \NN^n} { 1 \over \alpha !} (-\hbar)^{|\alpha|} 
\partial^{\alpha}_x D^{\alpha}_\xi a(z, \hbar \xi) e^{i (x-z )\cdot \xi }
d\xi d x d\eta
 $$
 $$
 = (2\pi)^n \int_{\RR^n} e^{i ( \varphi(y, \eta) -z\cdot \eta)}\,b(y,\eta)  \,
 \sum_{\alpha \in \NN^n} { 1 \over \alpha !} (-\hbar)^{|\alpha|} 
   \partial^{\alpha}_x D^{\alpha}_\xi  a(z, \hbar \eta) d\eta
 $$

 Now for $ \alpha \in \NN^n$ we set 
 $$
 c_\alpha(z; \hbar \eta ) = \partial^{\alpha}_x D^{\alpha}_\xi  a(z, \hbar \eta)
 $$ and we consider
 $$
 H_\alpha (y,z) = 
\int_{\RR^n} e^{i ( \varphi(y, \eta) -z \cdot \eta)}\,b(y, \eta) \,
c_\alpha(z, \hbar \eta ) d \eta
 $$ 
If we replace $ c_\alpha(z, \hbar \eta )$  by its Taylor expansion 
 $$
\sum_{\beta \in \NN^n} { 1 \over \beta !} \partial^{\beta}_z
c_\alpha(\varphi'_\eta ; \hbar \eta)
  (z-\varphi'_\eta)^{\beta}
 $$ 
then, using integrations by parts as above, it follows 
 easily that $H_\alpha (y,z) $
  is  the sum of a smooth function and of
  $$
\int_{\RR^n} e^{i ( \varphi(y, \eta) -z\cdot \eta)}\,\sum_{\beta \in \NN^n} {
1 \over \beta !} 
D^\beta_\eta \bigl(\, b(y, \eta) \, \partial^{\beta}_z
c_\alpha(\varphi'_\eta; \hbar \eta ) \,\bigr) d \eta
  $$ 
We observe  that if we apply the Leibniz rule for the term
$D^\beta_\eta \bigl(....\bigr) $ in the previous integral 
  then the following differential operators will appear
  \begin{equation} \label{diffop}
  D^{\beta - \gamma}_\eta b(y,\eta) \, D^{\gamma - \gamma'}_\eta
 (\varphi'_\eta)  D^{\gamma'}_\eta
   \partial^{\beta}_z
  \end{equation}
It is clear from Lemma \ref{lemma} that, expressed in the coordinates $(
\varphi'_\eta(y, \eta)\,,\, \eta )$, 
  these differential operators (\ref{diffop}) are 
$\E-$differential operators. Therefore we have just proved that $T(y,z)$
is the sum of a smooth 
  function and of:
  $$
  \int_{\RR^n} e^{i ( \varphi(y, \eta) -z\cdot \eta)}\, 
  \sum_{k\in \NN} \hbar^k P_k(a)(\varphi'_\eta(y, \eta) , \hbar\eta) d \eta
  $$ 
where the $P_k$ are $\E-$differential operators. 
  
Now we recall that the Schwartz kernel of $\Phi^*$ is the finite sum of
a smooth function and of terms 
of the type (\ref{Schwartz2}). So in order to analyze $\Phi \circ { \rm
Op}(a_\hbar) \circ\Phi^*$ it is 
enough to study the operator $K\circ { \rm Op}(a_\hbar) \circ K^*$ whose
Schwartz kernel is the finite 
  sum of a smooth function and of integrals of the type:
  $$
\int \int \int e^{i ( \varphi(y, \eta) -x\cdot \eta)} e^{-i ( \varphi(y',
\eta') -x\cdot \eta')}
  P_k(a)(\varphi'_\eta(y, \eta) , \hbar\eta) b_1(y',\eta') dx d\eta' d \eta 
  $$
  \begin{equation} \label{term}
  =(2\pi)^n \int_{\RR^n}  e^{i ( \varphi(y, \eta) -\varphi(y', \eta)}  
  P_k(a)(\varphi'_\eta(y, \eta) , \hbar\eta) b_1(y',\eta) d\eta
  \end{equation}
Moreover  we can write $\varphi(y, \eta) - \varphi(y', \eta) =
(y-y').\hat{\eta}(y,y',\eta) $ 
  where $\hat{\eta}(y,y,\eta) = \varphi'_y(y, \eta) $ 
and we can assume (at the expense of shrinking the local coordinates
charts) that 
$ \eta \rightarrow \hat{\eta} (y, y', \eta) $ is a local diffeomorphism
whose inverse 
is denoted $\hat{\eta} \rightarrow \eta (y, y',\hat{\eta} )$. With
these notations, we set:
  $$
A_k(y,y',\hbar, \hat{\eta}) = P_k(a)(\varphi'_\eta(y, \eta)\,,\,
\hbar\eta) b_1(y',\eta)
  $$ 
Then a change of variable formula allows us to see that the oscillatory
integral (\ref{term}) is equal to
  $$
(2\pi)^n \int_{\RR^n} e^{i (y -y')\cdot \hat{\eta} }\, A_k(y,y',\hbar,
\hat{\eta}) 
  \,|{D\eta \over D\hat{\eta}}| d\hat{\eta}
  $$ 
We observe that, expressed in the coordinates $(\varphi'_\eta(y,
\eta)\,,\, \eta)$, the vector fields 
  $ \partial_\eta (\varphi'_\eta) \partial_y  $ are $\E-$differential operators.
   Therefore one proves  easily the assertion of Part 2) of the Theorem 
  by replacing $A_k(y,y',\hbar, \hat{\eta}) $ by its Taylor expansion 
$$
\sum_{\beta \in \NN^n} {1 \over \beta !} \partial_{y'}^{\beta}
A_k(y,y,\hbar, \hat{\eta})_{|y=y'}\, (y'-y)^\beta 
   $$ 
and using, as before, integration by parts.

  \end{proof}

\section{The formal deformation and traces on $\cB X \cup_\phi \cB Y$ and
  regularized  traces on $\Psi$DO's}

   Recall first that $C^{+\infty}(M) $ is exactly the set of functions 
   $(f,g) \in C^{+\infty}( \cB X) \times C^{+\infty}( \cB Y) $ such that 
   $f - g \circ \phi' $ vanish of infinite order at the boundary of $\cB X$.

   We are going to use    the  $\star$-products denoted $\star_X$,  $\star_Y$ 
   on
   $\cB X$ and  $\cB Y$ defined in Propositions \ref{ext} and \ref{star5}. We
   set  
   $\Ah (\cB X ) = C^{+\infty}( \cB X) \otimes_\CC \CC[[\hbar]]$ and  
  $\Ah (\cB Y ) = C^{+\infty}( \cB Y) \otimes_\CC \CC[[\hbar]]$ .

   Let $\Ah ( M )$ be the vector space 
   given by
  $$
\{ (a,b) \in \Ah (\cB X ) \times \Ah ( \cB Y
   )  \mid  \tilde{\Phi} (\overline{a}) = \overline{b} \}
    $$ 
where $ \overline{a}$ 
    (resp. $\overline{b} $) denotes an element of 
    ${\mathbb B}^\hbar (\cB X) $ (resp. ${\mathbb B}^\hbar (\cB Y) $) 
    induced by $a $ (resp. $ b$).
    Theorem \ref{quantized} shows that $\Ah ( M )$ is 
   an algebra with respect to the diagonal product: $(\star_X,\star_Y)$. 
In particular, pairs of the form $(\,\sigma (\Phi^* \Phi),\, \sigma (\Phi
\Phi^*)\,)$ belong 
    to $\Ah ( M )$.

    In the statement of the next  proposition we will use the notations of
    Theorem \ref{quantized}. 

  \begin{proposition} \label{definition of deformation}
{\item 1)} Let $\chi \in C^{+\infty}(T^*X, [0,1])$ be such that $\chi(x,
\xi) = 0$ for 
$ || \xi|| \leq 1/2$ and $\chi(x, \xi) = 1$ for $ || \xi||\geq 1$. For
any $f \in C^{+\infty}( \cB X)$ we have :
  $$
\tilde{\Phi}( \chi f ) - ( \chi f )\circ \phi^{-1} \in \hbar {\mathbb
B}^\hbar (\cB Y )
  $$
  
{\item 2)} For each $b \in C^{+\infty}(\cB Y)$ one defines 
$b^-\in C^{+\infty}(\cB Y) $ by setting $b^-(\eta) = b(-\eta)$ 
for any $\eta \in \cB Y$.  The following formula:
  $$
  \forall (a,b)\in \Ah ( M ),\; \cal{{U}}(a,b) = 
(a + \sum_{k \geq 1}\, \hbar^k\, D_k(a), b^- ) \in C^{+\infty}(M)
\otimes_\CC \CC[[\hbar]]
  $$ defines a $\CC[[\hbar]]-$linear isomorphism $ \cal{{U}}$ from $\Ah ( M ) $
  to 
$C^{+\infty}(M)
\otimes_\CC \CC[[\hbar]] $.
{\item 3)} The product $\cal{{U}}(*_X,*_Y)$ defines an $\E-$deformation
of $M$ (or a star product) 
associated to  the symplectic Lie algebroid $(\E, [,], \omega)$ ( see
\cite{N-T 2} 
section 3.3).
  \end{proposition}
\begin{proof} Parts 1) and 2) are left to the reader. Part 3) is an easy
consequence of part 2) and of 
   Theorem \ref{quantized} 2).
\end{proof}

\begin{definition}\label{RR}
\item 
$\fA (M) $ denotes the formal deformation of $M$ associated to the symplectic
Lie algebroid $(\E, [,], \omega)$ constructed in the proposition
\ref{definition of deformation}.
\item 
The linear functional
$$ 
\tau_{\rm {can}} : \Ah ( M ) \rightarrow \CC[\hbar^{-1},\,\hbar]]
$$
is given by
   \begin{equation} \label{tau}
\forall (a,b) \in \Ah ( M ), \;\tau_{\rm {can}}( a, b) = \left\{
\begin{array}{c}
\mbox{ asymptotic expansion at $\hbar=0$ of }\\
 \hbar \mapsto \tau ( {\rm Op}(a_\hbar),\, {\rm Op}(b_\hbar) )
\end{array}
\right\}
   \end{equation}
where $\tau $ is the trace defined in Proposition \ref{regtrace}.
It follows immediately from the definition that $\tau_{\rm {can}}$ is a trace. 
\end{definition}

\noindent {\bf Computation of $\tau $.  }

  Since
the space of traces 
   on $ \Ah ( M )$ may be very big we introduce the following algebra:
   $$
   \DD^\hbar (M) = \Ah ( M ) \bigl[ \,\bigl(\,\, \chi\, ||\xi||,\,
( \chi \, ||\xi|| + \sum_{k \geq 1} \hbar^k D_k(\chi \, ||\xi|| ) )
\circ \phi^{-1}\,\bigr)\, \,\bigr] 
   $$

Another way of describing $  \DD^\hbar (M)$ is given by glueing from 

$$
{\tilde {\phi}}: \Ah (\cB X  )[\chi ||\xi||] \to  \Ah (\cB Y  )[\chi ||\xi||]
$$ 
 where $\chi$ is as in previous Proposition. It is easily seen that $\tau_{\rm
   {can}}$ defines,  by the same formula as (\ref{tau}), a trace on $ \DD^\hbar
 (M) $.

 Next Proposition 
   describes the space of   traces on  $\DD^\hbar (M) $.
   \begin{proposition} \label{trace} The space of traces with values in
   $\CC[\hbar^{-1},\,\hbar]]$  on 
   the algebra
    $ \DD^\hbar (M)$ is two dimensional over $\CC[\hbar^{-1},\,\hbar]]$. A basis 
    is given by $(\tau_{\rm {can}} , \tau_1 )$ where for any 
    $ ( a,b) \in \DD^\hbar (M)$ $\tau_1( a,b) = { \rm { Res}}_W (a) $. Here
     $ { \rm{ Res}}_W$ denotes Wodzicki's noncommutative residue.
   \end{proposition}
   \begin{proof} For $Z=X$ or $Y$ we set:
   $$
   C^{\infty}_0[\cB Z) \otimes_\CC \CC[[\hbar]]\,=\, \Ah_0(T^*Z)
   $$ 
where $ C^{\infty}_0[\cB Z)$ denotes the set of smooth functions which
vanish of infinite order 
   at $\cB Z \setminus T^*Z$. Then we have the following exact sequence:
   $$
   0 \rightarrow \Ah_0(T^*X) \oplus \Ah_0(T^*Y)\rightarrow \DD^\hbar
   (M)\rightarrow  {\mathcal T}(M) \rightarrow 0 
   $$ 
of $ \CC[[\hbar]]-$algebras. Here ${\mathcal T}(M)$
denotes the induced formal $\mathcal
E$-deformation of the sphere at infinity. A direct   construction of
this
 deformation may be described  as follows. Let ${\mathcal P}^i$ denote the space of
pseudodifferential
operators on, say, $X$ of order $\leq i$ modulo the
smoothing operators. Then the space of doubly infinite sequences
$$
\{P_i \}_{i\in \ZZ}, \ P_i \in {\mathcal P}^i ,\mbox{ there exists } i_0,
  P_i \in
{\mathcal
  P}^{i_0} \mbox{ for  $i$ large} 
$$
is a flat module over
$\CC[[\hbar]]$, where the multiplication by $\hbar$ acts as the right
translation. If we endow it with the product
$$
\{P_i \}_{i} \{Q_i \}_{i} =\{ \sum_{i+j=n} P_i Q_j \}_n 
$$ 
it is easily seen to be isomorphic to ${\mathcal T}(M)$.
 Any trace $\tau$ on $ {\mathcal T}(M)$ is
given by a sequence of $\CC$-linear, $\CC [[\hbar ]]$-valued functionals
$\tau_n$ on ${\mathcal P}^n$ such that
$$
\tau(\{ P_i\}) =\sum \tau_i  (P_i).
$$
The $\hbar$-linearity of $\tau$ implies that $\tau_{i+1} =\hbar \tau_i$ and
the trace condition on $\tau$ implies that each $\tau_i$ is a trace on the
algebra 
of pseudodifferential operators modulo the smoothing operators. Recall that, on
this latter algebra,  
the Wodzicki residue $res$ is 
  the unique trace up to
multiplicative constant. Thus $\tau$ is, up to
multiplicative constant, uniquely determined by $\tau_{-n} =res$, and hence
the space of traces on ${\mathcal T}(M)$ is one-dimensional.

We recall that $ \Ah_0(T^*X)$ is H-unital (in the sense of Wodzicki, see
\cite{Wo}), so we have 
   the following long exact sequence in cyclic cohomology:
   $$
   0 \rightarrow HC^0 ({\mathcal T}(M) )\rightarrow HC^0(\DD^\hbar (M) )
   \rightarrow 
   HC^0( \,\Ah_0(T^*X) \oplus \Ah_0(T^*Y)\,) \rightarrow
   $$
   $$
   HC^1({\mathcal T}(M) )\rightarrow \ldots
   $$ 
   From Section 2.1 we recall that the space of $\CC[[\hbar]]-$linear traces
   (with values in  
   $\CC[\hbar^{-1},\hbar]]$) on $ \Ah_0(T^*X)$ 
   is one dimensional and generated by ${\rm Tr}^X_{{\rm can}}$. By above,
   $HC^0 ({\mathcal T}(M) ) $ is one-dimensional. The connecting map:
   $$\delta:\, HC^0( \,\Ah_0(T^*X) \oplus \Ah_0(T^*Y)\,) \rightarrow
   HC^1({\mathcal T}(M)) 
   $$ 
is given by taking a trace on $\Ah_0(T^*X) \oplus \Ah_0(T^*Y)$, extending it
to a linear functional on $\DD^\hbar (M)$ and taking its Hochschild
boundary. In particular, it is not zero (this is equivalent to existence of a
pseudodifferential operator with nonzero index!). This implies that $
   HC^0(\DD^\hbar (M) )$ is either one or two dimensional. Since, with the 
   notations of the Proposition, $\tau_{{\rm can}}, \tau_1$ are two
   linearly independent  elements 
   of the vector space of traces on $\DD^\hbar(M)$,  
     the rest of the statement of above Proposition follows. 
   \end{proof}

\section{The algebraic index theorem for the Lie algebroid $\mathcal E$}  
 
The following Theorem is proved in \cite{N-T 2} and is an extension to
the symplectic Lie algebroid 
 $(\E,[,],\omega)$ of the 
 Riemann Roch theorem (on symplectic manifolds) for periodic cochains 
 of \cite{B-N-T}, \cite{N-T 1}.
 \begin{theorem} \label{RR} The following diagram is commutative:
 $$
\xymatrix{ 
  CC^{{\rm per}}_*(\Ah (M))
 \ar[r]^{\sigma} \ar[dr]_{\mu^\hbar } 
    & CC^{{\rm per}}_*(C^{\infty} (M)) \ar[d]^{\mu \cup 
 \hat{A}(M) \cup {\rm e}^\theta} \\
 &   (\,{}^{\E}{\Omega}^{2n-*}
 (M)[\hbar^{-1},\hbar]], d \,) & }
 $$ 
where $\sigma$ is the specialization map at $\hbar =0$,
  $\mu$ is the Hochschild-Kostant-Rosenberg map, $\mu^\hbar $ is the 
  trace density map 
defined in  \cite{N-T 1}
and  $\theta = {1 \over \sqrt{-1} \hbar} \omega 
  + \sum_{k\geq 0} \hbar^k \theta_k \in \uE{H^2}(M, \CC[[\hbar]])$ is the 
characteristic class of the deformation of the symplectic Lie algebroid
$(\E,[,],\omega)$ (\cite{N-T 2}). 
 \end{theorem}
 
The natural injection $ \Ah (M) \rightarrow \DD^\hbar (M)$ induces a
natural map:
 $$
CC^{{\rm per}}_*(\Ah (M)) \rightarrow CC^{{\rm per}}_*(\DD^\hbar (M))
 $$
Since the traces $\tau_{\rm {can}} $ and $\tau_1 $ of Proposition \ref{trace},
and the 
trace density map $ \mu^\hbar$ extend  to $CC^{{\rm
per}}_*(\DD^\hbar (M)) $, they can be identified using the following
result.  
 
 \begin{proposition} \label{equality} 
{\item 1)} The $\CC-$vector space ${}^{\E}{H}^{2n}( M, \CC)$ is
 two-dimensional. The vector space of $\CC-$linear forms on
 ${}^{\E}{\Omega}^{2n}(M)$ which vanish on the range of the 
 $\E-$exterior derivative  ${}^{\E}{d}$
admits a unique linear  basis  $({}^{{\rm reg}}{\int}, \int_1 )$, 
characterized by the following properties

For any $(\alpha, \beta) 
\in {}^{\E}{\Omega}^{2n}(M)$ such that $\alpha$ (resp. $ \beta$) is zero
in a neighborhood of the 
 boundary of $\cB X$ (resp. $ \cB Y$) 
 $$
{}^{{\rm reg}}{\int}(\alpha,\beta ) \, =\, \int_{T^* X} \alpha -
\int_{T^* Y} \beta, 
 \; \;\int_1 (\alpha,\beta ) = 0 .
 $$ 
Moreover $\int_1 \circ \mu^\hbar = { \rm{ Res}}_W$ .
 {\item 2)} There exists a  constant $C$ such that
 $$
\tau_{\rm {can}} = {}^{{\rm reg}}{\int} \circ \mu^\hbar + C \int_1 \circ
\mu^\hbar
 $$ 
Moreover $\int_1 \circ \mu^\hbar(1,1) =0$ and for any $(a,b) \in
\DD^\hbar (M)$ such that 
$a$ is zero in a neighborhood of the boundary of $\cB X$, $\int_1 \circ
\mu^\hbar(a,b) =0$.
 
   \end{proposition}
\begin{proof} 

\noindent 1) A standard Mayer-Vietoris sequence argument shows that
${}^{\E}{H}^{2n}( M, \CC)$  
is indeed two dimensional. The fact that $({}^{{\rm reg}}{\int}, \int_1
)$ defines a basis is left to the reader. 

\noindent 2) This is an easy consequence from part 1) and of the properties
(see \cite{N-T 1}, \cite{N-T 2} ) of the trace density map $\mu^\hbar$. 
 \end{proof}

\section{Local formula for the index of a Fourier Integral Operator}

\begin{theorem} 
\item
Let $\Phi$ be a Fourier Integral Operator and $\fA (M)$ the formal deformation
of $M$ associated to it as in Definition \ref{RR}. Then
$$
\ind \Phi\,=\, \int_M {\rm e}^{\theta_0} \hat{A}(M) ,
$$
where $\theta_0$ denotes 
the characteristic class of the deformation of the 
 Lie Algebroid (${\mathcal E}, [\,,\,], \omega $) given by $\fA (M)$. 
\item
Let $\nabla_X$ be a connection $\nabla_X$ on the tangent bundle $T( \cB X)$
and $\hat{A}(T^*X)$ an associated representative form of  the
$\hat{A}-$class of $\nabla_X$.
The symplectomorphism 
$\phi$ induces a connection $\phi_* (\nabla_X)$ on the tangent space of 
${\cB} Y \setminus B^*(Y)$. 
Let $\nabla_Y$ denote its  extension to a connection of $T({\cB Y})$
and $\hat{A}(T^*Y)$  an associated representative differential form  of the
$\hat{A}$-class of $\nabla_Y$. 
Then  
$$
 \ind \Phi\,=\, \int_{B^*(X)} {\rm e}^{\theta_0} \hat{A}(T^*X) \,-\,
\int_{B^*(Y)} {\rm e}^{\theta_0} \hat{A}(T^*Y) 
$$
   \end{theorem} 
 \begin{proof} 
 One obtain this formula by first applying
    Proposition \ref{regtrace},
   Theorem \ref{RR} and Proposition \ref{equality} and 
   then by letting $\hbar \rightarrow 0^+$. As we will see below, the involved
    characteristic classes  of vector bundles on $M$ are in fact
    standard de Rham cohomology classes and hence  the regularized integral
    coincides with the orientation class of $M$.
\end{proof} 
   The previous formula shows that if $\phi$ extends as a symplectomorphism 
   $T^*X \rightarrow T^*Y$ up to the zero section then $\,\ind \, \Phi =\, 0$.

For a deformation associated with a
Fourier Integral Operator (as in  Proposition \ref{definition of
  deformation}) the characteristic
class $\theta$ of Theorem \ref{RR} is in fact of the form
$$
\theta = {1 \over \sqrt{-1} \hbar} \omega + \theta_0
$$
where $\theta_0 \in H^2(M, \CC)$ is a closed  differential form (not only an
${\E}-$differential form).
 In order to do this and to identify the relevant characteristic class 
 we will give below a slightly nonstandard description of a formal
 deformation.

\subsection{General construction of 
the characteristic class of a formal deformation}

Let us start with some notation.

Let ${\mathbb A}^ \hbar $  denote the Weyl algebra 
of the symplectic vector space
$\RR^{2n}$ with the  standard symplectic structure, i.e. the algebra  generated by
the vectors $\fx_l ,\fxi_l $ ($1 \leq l \leq n$) satisfying the  relations  $[\fxi_k ,\fx_l] = \sqrt{-1}
 \hbar
\delta_{k,l}$.  The algebra ${\mathbb A}^ \hbar  $ is completed in the topology associated
to the ideal generated by $\{\fx_l ,\fxi_l ,\hbar ;\ 
1 \leq l \leq n\}$
 and has the grading induced by  
$$
deg \,\fx_l =deg\,\fxi_l =1,\ deg \,\hbar =2 .
$$
The corresponding Lie algebra $\frac{1}{\hbar }{\mathbb A}^ \hbar $ will be
denoted by $\tilde{\frak g}$. We set: 
$$
{\frak g}=\mbox{Der }({\mathbb A}^ \hbar ) =\tilde{\frak g}/\mbox{center}
$$ 
  and
 $$
G = \mbox{Aut}({\mathbb A}^ \hbar )= \exp ({\frak g}_{\geq 0})
$$
We set
$$
\tilde{G} =\{ g \in \frac{1}{\hbar }{\mathbb A}^ \hbar 
\mid g\in {\frak sp}(2n, \RR ) \ mod \,{\frak g}_{\geq 1} \}  
$$
and will endow it with the group structure coming from the exponential
map. Note that $\tilde G$ is an extension of $G$ associated to the (Lie
algebra) central extension $\tilde{\frak g}$ of  $\frak g$.
 
We endow the bundle $\RR^{2n} \times \fA$ with the obvious fiber-wise
action of
$\tilde{G}$ and with 
 the $\tilde g$-valued (Fedosov) connection  
$$
\tilde{\nabla}^0=
\sum_{l=1}^{n}(\,
 d\xi_l (\partial_{\xi_l} -{ 1\over \sqrt{-1} \hbar }\fx _l )+ 
 dx_l (\partial_{x_l} + { 1\over  \sqrt{-1}\hbar}\fxi_l )\,) .
$$
Let us recall (see Section 2.2) that a local chart of 
$\cB \RR^n$ near $ \cB \RR^n \setminus T^* \RR^n$ is given by:
\begin{equation} \label{chartbis}
(\,x_1,\ldots, x_n; t= {1\over ||\xi||}, \theta=(\theta_1, \ldots 
\theta_{n-1})\,)\, t\geq 0,\; \theta \in S^{n-1}
\end{equation}
By using the local coordinates (\ref{chartbis}) one checks easily 
that $\tilde{\nabla}^0 $ extends as an ${{\E}^{\RR^n}}-$connection,
still denoted $\tilde{\nabla}^0 $,
of $\cB \RR^n\times \fA$.  

\bigskip

The description given below of a formal deformation of a symplectic
Lie algebroid structure on $M$ is just the representation of the
Fedosov construction in terms of the bundle of jets on $M$ with the
fiber-wise product structure induced by the *-product (which is
isomorphic to Weyl bundle) .

\bigskip

\noindent {\bf Local description of the characteristic  class $\theta$ of a
  formal deformation.}

\bigskip

\begin{em}
 The
deformation is described by a local (Darboux) 
 cover  $\{U_i \}_{i \in I}$ of $(M,\omega )$,
  a  collection of functions
 $\{ g_{i,j}:U_i \cap U_j \rightarrow  \tilde{G}
\}$  and a collection of $\tilde{\frak g}$-valued  ${{\E}}-$connections
$\tilde{\nabla}_i $ on  $U_i \times \fA$ which, when expressed in
terms of local Darboux coordinates
$(x_1 ,\ldots ,x_n ,\xi_1 , \ldots
,\xi_n )$ (resp (\ref{chartbis}))  if $ U_i$ does not meet (resp meets)
the boundary at infinity 
, are equal to $\tilde{\nabla}^0$
modulo $\tilde{\frak g}_{\geq 1}$ and so that the three
 following conditions
hold.

1) The cocycle condition: 
$$ g_{i,j}g_{j,i} = 1\; {\rm and}\; g_{i,j} g_{j,k} = g_{i,k} 
\; {\rm on}\; U_i \cap U_j \cap U_k
$$
In particular $\{ g_{i,j}:U_i \cap U_j \rightarrow  \tilde{G}
\}$ define a smooth bundle $\mathcal{W}$ of algebras over $M$ with 
fiber isomorphic to ${\mathbb A}^ \hbar $ and the structure group $\tilde{G} $.

2) The local connections $\tilde{\nabla}_i$ define a $\tilde{\frak
  g}$-valued connection 
$\tilde{\nabla}$ on the bundle $\mathcal{W}$, i.e.:
$$
g_{i,j} \tilde{\nabla}_j = \tilde{\nabla}_i g_{i,j} 
$$

3) The induced ${\mathfrak  g}$-valued connection ${\nabla}$ on the bundle
$\mathcal{W}$ is flat, i. e. $\theta =\tilde{\nabla}^2$ 
is a globally defined differential form on $M$ with values in the
center of $\tilde{\frak g}$, necessarily of the form
$$
\frac{1}{ \sqrt{-1} \hbar} \omega +\theta_0  \ \mbox{ where } 
\theta_0 \in  \Omega^2 (M, \CC[[ \hbar]] ). 
$$

The algebra of $\nabla$-flat sections of $\mathcal {W}$ is a formal 
deformation of $(M,\omega)$ whose characteristic class is $\theta$.

\end{em}

%%%%%%%%%%%%%%%%%%%%%%%%%%%%%%%%%%%%
%%%%%%  construction of deformation
%%%%%%%%%%%%%%%%%%%%%%%%%%%%%%%%%%%%

\subsection{Local canonical liftings}

\noindent We endow  $\RR^{2n}$ with its canonical symplectic structure
$\omega= \sum_{l=1}^{n} d\xi_l \wedge dx_l$. 
Given any smooth, $\CC [[ \hbar ]]$-valued function $H$ on $\RR^{2n}$, we set
$$
H_0 = H(x,\xi )|_{\hbar =0},\ H_1 = \sum_{l=1}^n (\fx_l \partial_{x_l}H_0
 +\fxi_l
\partial_{\xi_l}H_0 )    
$$
and
$$
\tilde{H}=\sum \frac{\fx^\alpha \fxi^\beta }{\alpha ! \beta
  !}\partial_{x}^\alpha \partial_{\xi}^\beta H  .  
$$

We will associate to $H$ the following $\tilde{\frak g}$-lift of the Lie
derivative  
${\mathcal L}_{\{ H,\ \}}$:
$$
{\mathcal D}_H = {\mathcal L}_{\{ H,\ \}} + \frac{1}{\hbar} (\,
\tilde{H}-H_{0} - H_1 +\frac{1}{2} \hbar \sum_{l=1}^n
 \partial^2_{x_l ,\xi_l} H_0
\, ) . 
$$
We can think of it as an element of the Lie algebra of the semidirect product
of $C^\infty  (\RR^{2n} ,\tilde{G})$ by the pseudogroup of local diffeomorphisms
of $\RR^{2n}$. The $\mathcal D_H$'s form a Lie algebra, in fact  
$$
[{\mathcal D}_H  ,{\mathcal D}_K ]= {\mathcal D}_{\frac{1}{\hbar} ( H * K - K
  * H )} , 
$$
and they satisfy
$$
[{\mathcal D}_H  , \tilde{\nabla}^0 ]= 
-\frac{1}{2} d ( \sum_{l=1}^n \partial^2_{x_l
  ,\xi_l} H_0  ). 
$$
We will also have an occasion to use 
\begin{equation} \label{use}
{\mathcal D}_H^0 = {\mathcal L}_{\{ H,\ \}} + \frac{1}{\hbar} (\,
\tilde{H}-H_{0} - H_1  \,) , 
\end{equation}
which commutes with $\tilde{\nabla}^0$.

\bigskip

\subsection{The cotangent bundle case}

The deformation of  $T^*X$ associated to the sheaf of differential operators
on $X$ can be now described as follows. 

Locally on a coordinate domain
$U\subset X$ we use coordinates on $U$ to give an explicit
symplectomorphism 
$$
T^* U \rightarrow \RR^{2n}
$$
and use Weyl deformation of $\RR^{2n}$ to construct the deformation of
$T^*U$. This amounts  to  the choice of a ($\tilde{\frak g}$-valued) connection
given in our local coordinates $(x_1 ,\ldots ,x_n )$ on $U$ and the
induced local coordinates $(x_i ,\xi_i)_{i=1,...n}$ on $T^*U \simeq
U\times \RR^{n}$ by  
$$
\tilde{\nabla}^0 = d+ \frac{1}{\hbar}\sum_{l=1}^n 
( dx_l  \fxi_l - d\xi_l \fx_l ).
$$
The infinitesimal change of coordinates on $U$ is given by a vector field of
the form $\sum_{l=1}^n X_l \partial_{x_l}$ and the associated infinitesimal
symplectomorphism of $T^* U$ is given by the Hamiltonian vector field $\{
\sum_{l=1}^n X_l \xi_l , \cdot\}$.  

It is immediate to see that the map
$$
\sum_l X_l \partial_{x_l} \mapsto  {\mathcal D}_{\sum_l X_l \xi_l}
$$
is the Lie algebra homomorphism.

The associated local diffeomorphisms (coordinate changes)
$\exp \sum_l X_l \partial_{x_l}$ lift to a local isomorphisms of the bundle
$T^* U\times  {\mathbb A}^ \hbar$ given by 
$\exp {\mathcal D}_{\sum_l X_l \xi_l}$.

Given a local coordinate cover $\{  U_i \}_{i\in I}$ of $X$ it is now immediate to
construct the associated $\tilde{G}$-valued cocycle $\{ g_{ij} \}$ glueing the
bundles together. Note that, since $\mathcal D$'s do not commute with the
connection $\tilde{\nabla}^0 $, the corresponding collection of connections 
$$
\tilde{\nabla}_i = \tilde{\nabla}^0 \ \mbox{in i'th coordinate system
  on $T^*U_i$ }  
$$
do not glue together. But it is not difficult to check that 
$$
g_{ij} \tilde{\nabla}_i g_{ji} = \frac{1}{2} d \log \det  Dg_{ij} ,
$$
where $Dg_{ij}$ is the induced action  of $g_{ij}$  on the tangent
bundle.  By trivializing the cocycle 
$\frac{1}{2} d \log \det  Dg_{ij}$ in $\check {C}^1 (X, \Omega^1 (T^* X))$ we
get a globally defined connection $\tilde{\nabla}$ and it is immediate that
the characteristic class of the associated deformation is $\frac{1}{2}
\pi^* c_1 (T_\CC X)$,  where $\pi:T^*X \rightarrow X$ is the canonical
projection. 
It is also immediate that the deformation constructed in this way coincides
with the one associated to the calculus of differential operators on $X$,
while its jet at $\xi =\infty$ gives the deformation associated to the
calculus of pseudodifferential operators on $X$, the characteristic
class being given by the jets at $\xi = \infty$ of
 ${1\over 2}c_1(\E_\CC)$ (recall that the real symplectic vector
 bundle  $\E$ is the realification of 
a complex vector bundle $\E_\CC$).   
\bigskip

\subsection{The Lie algebroid case}

Recall now that the  Lie algebroid  (on $M$) $(\E, [\,, \,], \omega )$ 
 is  given by glueing (at infinity) 
 the two cotangent bundles $(T^* X,
 \omega^X)$
 and $(T^* Y, \omega^Y)$
 by
the  symplectomorphism $\phi'$. 
 To
construct the deformation in this case, we will use the following
 data, whose existence follows immediately from the compactness of the
 co-sphere bundles of $X$ and $Y$.  
\begin{enumerate}
\item A  local coordinate cover $\{ U_i \}_{i \in I}$ of $X$ 
and an open relatively compact
  neighborhood $U_X$ of the zero section in $T^* X$; 
\item A  local coordinate cover $\{ V_i \}_{i \in I}$ of $Y$
and an open relatively compact neighborhood $U_Y$ of the zero section in
 $T^*Y$; 
\item For each $i \in I$ a one-homogeneous real-valued function $H_i$ on
  $T^*X \setminus X \cong   T^*Y \setminus Y$ such that the restriction
  $\phi_i$ of the symplectomorphism $\phi$ to $T^* U_i \setminus U_X$
  is  given by  integrating the (time dependent) hamiltonian flow
  ${\mathcal L}_{H_i}$.   
 \end{enumerate}

Using the above data, we can construct cocycles
$$
T^* U_i \cap T^* U_j \ni p \rightarrow g_{ij} (p)\in C^\infty (T^* U_i \cap
T^* U_j  ,\tilde{G}) 
$$
and
$$
T^* V_i \cap T^* V_j \ni p \rightarrow h_{ij} (p)\in C^\infty (T^* V_i \cap
T^* V_j  ,\tilde{G}) 
$$
intertwining the flat connections
$\tilde{\nabla}^0_i$   up to the term $\frac{1}{2}d \log \det  Dg_{ij}$
($\frac{1}{2} d \log \det  Dh_{ij} $ respectively) as in the cotangent bundle
case. We can also construct, using the notation of (\ref{use}),
 a lifting of  $\phi_i$ 
$$
\exp {\mathcal D}^0_{H_i} = \Psi_i : \Gamma (T^*U_i\setminus U_X
 , \fA )\rightarrow \Gamma
(T^*V_i \setminus U_Y , \fA ) .
$$
>From now on  we view   $\Phi_i$ as local isomorphisms of jets at infinity of
the $\tilde{G}$-bundles on compactified cotangent bundles of $X$ and
$Y$ constructed from the cocycles $g_{i,j}$ and $h_{i,j}$.   
While both $g_{i,j}$'s and $h_{i,j}$'s do
 satisfy the cocycle conditions on $T^* (X)$ ($T^*(Y)$ respectively), however 
$$
\lambda_{ij} =\Psi_j^{-1} h_{ji} \Psi_i g_{ij} \neq 1 .
$$
and hence we do not yet have the data necessary to construct the
bundle $\mathcal W$ over $M$.

The following facts are easy corollaries of the construction: 
\begin{enumerate}
\item 
$\lambda_{ij} =1 \ \mbox{ mod } \tilde{G}_{\geq 1} $;
\item 
$\lambda_{ij}\tilde{\nabla}_j^0\lambda_{j,i}= \tilde{\nabla}_i^0 - \frac{1}{2}d
\log \det  Dg_{ij} -\frac{1}{2} d \log \det  Dh_{ij}$; 
\item 
$\lambda_{ij}$ form a two-cocycle with values in $\tilde{G}$.
\end{enumerate}

To begin with, note that both $\frac{1}{2}d \log \det  Dg_{ij}$ and
$\frac{1}{2}d \log \det  Dh_{ij}$ as cohomology classes on $T^*X\setminus X$
and $T^*Y\setminus Y$ represent (under our symplectomorphism) the same
cohomology class, to wit half of the first Chern class of the tangent bundle
with the complex structure induced by the symplectic form. Since these
vanish, we can find a 
zero-$\check{C}$ech cochain $\tau_i$ of the sheaf of   functions with
values in $\CC \setminus \{ 0\}$ and such that  $\tau_i \lambda_{ij}
\tau_j^{-1}$ intertwines the (local) flat connections
$\tilde{\nabla}_i^0$ and $\tilde{\nabla}_j^0$.  

In particular, $\tau_i \lambda_{ij} \tau_j^{-1}$ are given by exponentials of
jets of $\tilde{\nabla}^0_i$-flat sections of the bundle $T^*U_i \times \fA$
and, using a partition of unity, they can be written in the form 
$\tau_i \lambda_{ij} \tau_j^{-1} = \lambda_i
\lambda_j^{-1}$, where $\lambda_i$ 
is a jet of a flat section of the Weyl bundle supported on  $T^*U_i
\setminus U_X$.
We now define 
an operator $\Psi$, acting on the set of sections of the Weyl bundle 
 $\mathcal{W}$, by
setting 
for each $i \in I$: 
$$
\Psi_{|T^*U_i \setminus U_ X} =  \Psi_i \tau_i {\mathcal M}_{\lambda_i} .
$$
Here ${\mathcal M}_{\lambda_i}$ stands for the operator of multiplication with
the flat section $\lambda_i$.
It is easy to see that $\Psi$ descends to an isomorphism of jets at the
sphere at infinity of the deformations of cotangent bundles constructed above so
that
$$ \forall a \in C^\infty(\cB X),\;
\Psi (a)  = (\, a+ \sum_{k\geq 1} \hbar^k \hat{D}_k(a)\,) \circ \phi^{-1}
$$
holds in  ${\mathbb {B}}^\hbar(\cB Y)$, where the $\hat{D}_k$ are
$\E-$differential operators. 
Hence, as in Proposition \ref{definition of deformation}, $\Psi$ 
 induces a deformation of the  Lie algebroid $(\E, [\,,\,], \omega)$. 

\bigskip

\subsection{The characteristic class of $\fA (M)$}

The characteristic class of the deformation constructed above can now be
easily obtained as follows. The collection $g_{ij}, h_{ij}$, and the jet at
infinity of $ \Psi_i {\mathcal M}_{\lambda_i}$ give a cocycle with
values in $\tilde{G}$, and 
it commutes with local flat connections up to 
the Cech cocycle  given by the collection of differential forms
\begin{equation}\label{class}
\frac{1}{2} d\log \det Dg_{ij} , d \log \tau_i , \frac{1}{2} d\log \det
Dh_{ij}  . 
\end{equation}
 As in the case of cotangent bundle, we can correct
local connections by a scalar term. The characteristic class $\theta_0$ of the
deformation is given by (\ref{class}) as a cochain in  $\check{C}^1(M,
\Omega^1(T^*M) )$
 Moreover, in the case that both $X$ 
and $Y$
admit metalinear structures, the collection $\{ \tau_i\}_{i \in I}$
 can be thought of
as glueing of the pulled back)of the half-top form bundles of $X$ and
$Y$ along the graph of the symplectomorphism into a line bundle
$\mathcal L$ over $M$ and, in this case, 
$$
\theta_0 = c_1( \mathcal {L})
$$

\bigskip

\subsection{The Fourier Integral Operator}

To get the Fourier integral  operator we will work locally. We will
dispense with the half-density 
bundles (trivial in any case) for the sake of simplicity of
notation. We will begin by  constructing, for each i,   an operator on
$L^2(\RR^n )$ as follows.  
Choosing local coordinates on $U_i$ and $V_i$, we can assume that $H_i$
(introduced at the end of section 7.4) is
actually a smooth function on $T^* \RR^{2n}$ which is 1-homogeneous in the
cotangent direction.  The solution of the following  differential equation 
$$
\frac{d}{dt} T_i (t) = \mbox{Op}(\sqrt{-1}H_i) \circ  T_i (t), \ T_i (0)=1
$$ 
is a smooth family of bounded operators. Using the fact
that  $\{\tau_i\}_{i\in I}$ is a 
$\check{C}$ech zero- cochain of the sheaf of   functions with
values in the unit circle and  proceeding as in 
\cite{Fedosov}, one checks that  $T_i(1)$ 
satisfies 
$$
\mbox{Ad }(T_i (1) \mbox{Op}(\lambda_i ) ) \mbox{Op}(f_\hbar ) \sim \
\mbox{Op}( \Psi (f )_\hbar  )  
$$
mod ${\hbar}^{+\infty}$ as $\hbar \rightarrow 0$ whenever 
 supp$f\subset T^*U_i \setminus X$ ( recall that $\lambda_i$
 is introduced in Section
7.4). 
 In other words, the deformation constructed above is associated (in the 
 sense of Proposition \ref{definition of deformation})  to the almost 
unitary 
Fourier Integral Operator
 $\Phi = \sum_{i \in I}  T_i (1) \mbox{Op}(\lambda_i )$ whose canonical relation is $C_\phi$.
 Moreover, the index of this operator $\Phi_0$ is given by  
$$
\int_M \hat{A}(M)e^{\theta_0} ,
$$

\begin{remark}
The result above depends on the choice of the $\tau_i$'s
  which in turn determine
  the homotopy class of the
symbol of the Fourier Integral Operator. Moreover, since the characteristic
classes involved are given by differential forms associated to connections on
a vector bundle over $M$ and $\Omega \subset {}^\E \Omega$, the $\E$-classes
involved in the index formulas are in fact identical with corresponding
standard characteristic classes. 
\end{remark} 

Let us recall that 
the
real vector bundle $\E \simeq TM$ is given by realification of a complex
vector bundle $ \E_\CC$ on $M$ (the almost complex structure coming from the
symplectic vector bundle structure on $\E$). Moreover, it is easy to see, 
 that there exists a choice of the $\tau_i$'s such that the associated 
 characteristic
 class
of the deformation coincides with ${1 \over
    2}c_1(\E_\CC) $. This gives the  
 following result (compare with \cite{E-M}
and
\cite{W 2}):

\begin{theorem} There exists an almost unitary Fourier
  integral 
operator $\Phi_0$ (as in Section 3) whose canonical relation is
$C_\phi$ and such that:
$$
\ind \Phi_0 \,=\, \int_M \hat{A}(M)e^{{1 \over
    2}c_1(\E_\CC)}
$$
\end{theorem}

%%%%%%%%%%%%%%%%%%%%%%%%%%%%%%%%%%%%%%%%%%%%%%%%%%%%%%%%%%%%%%%%%%%%%%%%%%%%%%%%%%%%%%%%%%%%%%%%%%%%%%%%%%%%%%%%%%%%%%%%%%%%%%%%%%%%%%%%%%%%%%%%%%%%%%%%%%%%%%%%%%%%%%%%%%%%%%%%%%%%%%%%%%%%%%%%%%%%%%%%%%%%%%%%%%%%%%%%%%%%%%%%%%%
%%%%%%%%%%%%%%%%%%%%%%%%%%%%%%%%%%%%%%%%%%%%%%%%%%%%%%%%%%%%%%%%%%%%%%%%%%%%%%%%%%%%%%%%%%%%%%%%%%%%%%%%%%%%%%%%%%%%%%%%%%%%%%%%%%%%%%%%%%%%%%%%%%%%%%%%%%%%%%%%%%%%%%%%%%%%%%%%%%%%%%%%%%%%%%%%%%%%%%%%%%%%%%%%%%%%%%%%%%%%%%%%%%%
%%%%%%%%%%%%%%%%%%%%%%%%%%%%%%%%%%%%%%%%%%%%%%%%%%%%%%%%%%%%%%%%%%%%%%%%%%%%%%%%%%%%%%%%%%%%%%%%%%%%%%%%%%%%%%%%%%%%%%%%%%%%%%%%%%%%%%%%%%%%%%%%%%%%%%%%%%%%%%%%%%%%%%%%%%%%%%%%%%%%%%%%%%%%%%%%%%%%%%%%%%%%%%%%%%%%%%%%%%%%%%%%%%%
%%%%%%%%%%%%%%%%%%%%%%%%%%%%%%%%%%%%%%%%%%%%%%%%%%%%%%%%%%%%%%%%%%%%%%%%%%%%%%%%%%%%%%%%%%%%%%%%%%%%%%%%%%%%%%%%%%%%%%%%%%%%%%%%%%%%%%%%%%%%%%%%%%%%%%%%%%%%%%%%%%%%%%%%%%%%%%%%%%%%%%%%%%%%%%%%%%%%%%%%%%%%%%%%%%%%%%%%%%%%%%%%%%%
%%%%%%%%%%%%%%%%%%%%%%%%%%%%%%%%%%%%%%%%%%%%%%%%%%%%%%%%%%%%%%%%%%%%%%%%%%%%%%%%%%%%%%%%%%%%%%%%%%%%%%%%%%%%%%%%%%%%%%%%%%%%%%%%%%%%%%%%%%%%%%%%%%%%%%%%%%%%%%%%%%%%%%%%%%%%%%%%%%%%%%%%%%%%%%%%%%%%%%%%%%%%%%%%%%%%%%%%%%%%%%%%%%%

\end{document}